\documentclass[12pt]{amsart}
\usepackage{graphicx,amscd}
\input xy
\xyoption{all}

\title{Hilbert Series of Subspace Arrangements}
\author{Harm Derksen}
\thanks{The author was partially supported by NSF grant, DMS 0349019}
\newtheorem{theorem}{Theorem}[section]

\newtheorem{proposition}[theorem]{Proposition}

\newtheorem{corollary}[theorem]{Corollary}
\theoremstyle{definition}
\newtheorem{definition}[theorem]{Definition}
\newtheorem{example}[theorem]{Example}
\newcommand{\PP}{{\mathbb P}}
\newcommand{\C}{{\mathbb C}}
\newcommand{\R}{{\mathbb R}}
\newcommand{\Z}{{\mathbb Z}}

\newcommand{\HH}{{\mathcal H}}
\newcommand{\N}{{\mathbb N}}

\newcommand{\cd}{\operatorname{cd}}
\newcommand{\reg}{\operatorname{reg}}

\newcommand{\im}{\operatorname{im}}
\begin{document}
\maketitle
\begin{abstract}
The vanishing ideal $I$ of a subspace arrangement $V_1\cup V_2\cup \cdots \cup
V_m\subseteq V$
is an intersection $I_1\cap I_2\cap \cdots \cap I_m$
of linear ideals.
We give a formula for the Hilbert polynomial of $I$ if the subspaces
meet transversally.
We also give a
formula for the Hilbert series of the product ideal
$J=I_1I_2\cdots I_m$ without any assumptions on the subspace arrangement.
It turns out that the Hilbert series of $J$ is a combinatorial invariant
of the subspace arrangement: it only depends on the intersection
lattice and the dimension function. The graded Betti numbers
of $J$ are determined by the Hilbert series, so they
are combinatorial invariants as well.
We will also apply our results to
Generalized Principal Component Analysis (GPCA),
a tool that is useful for  computer vision and image processing.
\end{abstract}
\section{Introduction}
Suppose that $V$ is an $n$-dimensional $K$-vector space. A subspace
arrangement is a union
$$
{\mathcal A}=V_1\cup \cdots \cup V_m
$$
where $V_i$ is a subspace of $V$ for all $i$. 
Interestingly, various algebraic and topological properties of
the arrangement ${\mathcal A}$ only depend on the dimensions
$n_S:=\dim_K\bigcap_{i\in S} V_i, S\subseteq \{1,2,\dots,m\}$. 
Such properties are called {\it
combinatorial invariants\/} of the subspace arrangement.
 For example, if $K=\R$,
then the topological Betti numbers of the complement $V\setminus {\mathcal A}$
are combinatorial invariants (see~\cite{GM}).
If $K=\C$, then the cohomology ring of $V\setminus {\mathcal A}$
is a combinatorial invariant (see~\cite{LS,DGM}).
For more on subspace arrangements and hyperplane arrangements, see~\cite{OT}.


Let $I_j\subseteq K[V]$ be the vanishing ideal $V_j\subseteq V$
for $j=1,2,\dots,m$. The vanishing ideal $I$ of ${\mathcal A}$
is equal to the intersection $I_1\cap I_2\cap \cdots \cap I_m$.
We also define $J=I_1I_2\cdots I_m$ as the product ideal.
We give a formula for the
Hilbert series of $J$ (Theorem~\ref{thmmain1}).
 We also will give a formula
for the Hilbert polynomial of $I$ if all subspaces meet transversally
(Theorem~\ref{theo1}).
The Hilbert series of $J$ is a combinatorial
invariant (Corollary~\ref{cor5}), but the Hilbert polynomial of $I$ is not (Example~\ref{exlast}).
The Betti numbers and graded Betti numbers of $J$ are also combinatorial invariants
(Corollary~\ref{corComb}).

The original motivation for this paper comes from computer vision.
A generalization of Principal Component Analysis naturally
leads to the question of recovering the dimensions
$n_i:=\dim V_i$, $i=1,2,\dots,m$, given the Hilbert polynomial
of the subspace arrangement. For more on Generalized
Principal Component Analysis, see~\cite{VMS}.

{\bf Acknowledgement.}
I would like to thank Robert Fossum for interesting discussions,
for explaining me the results in computer vision, and for suggesting
to study the Hilbert function of a general subspace arrangement.
I also thank Allen Yang and Yi Ma for useful references.

\section{Hilbert functions, series and polynomials}
Suppose that $V$ is an $n$-dimensional vector space over a field $K$.
We identify the coordinate ring $R:=K[V]$
with the polynomial ring $K^{[n]}:=K[x_1,x_2,\dots,x_n]$ in $n$ variables
by choosing a basis in $V$.
 There is
a natural grading $R=\bigoplus_{d\in \N} R_d$ where
$\N=\{0,1,2,\dots\}$ is the set of natural numbers and
$R_d$ denotes the space of homogeneous polynomials of degree $d$. 
Let $\Z$ be the integers and suppose that $M=\bigoplus_{d\in \Z} M_d$ is a finitely
generated graded $R$-module. We have
$M_d=0$ for $d\ll 0$ because $M$ is finitely generated. 
The {\it Hilbert function\/} $h_M$ of $M$ is
$$
h_M(d)=h(M,d)=\dim_K M_d,\quad d\in \Z.
$$
The {\it Hilbert series\/} of $M$ is defined by
$$
\HH(M,t):=\sum_{d\in\Z}h(M,d)t^d.
$$
It is a Laurent series because $h(M,d)=0$ for $d\ll 0$.
\begin{example}
For $M=R$ we get
$h(R,d)={n-1+d\choose n-1}$ for $d\geq 0$ and $h(R,d)=0$ for $d<0$.
So we have
$$
\HH(R,t)=\sum_{d=0}^\infty {\textstyle {n-1+d\choose n-1}}t^d=\frac{1}{(1-t)^n}.
$$
Define the polynomial $\widetilde{h}_R$ by
$$
\widetilde{h}_R(d)=\widetilde{h}(R,d)=\frac{(d+n-1)(d+n-2)\cdots (d+1)}{(n-1)!}.
$$
One easily checks that $h_R(d)=\widetilde{h}_R(d)$ for $d\geq 1-n$.
\end{example}
Let $M$ be again a finitely generated graded $R$-module.
For $r\in \Z$ we
define the shifted module $M[r]=\bigoplus_{d\in \Z} M[r]_d$ by
$M[r]_d:=M_{r+d}$, $d\in \Z$.
Shifting the degrees affects the Hilbert function and the Hilbert series
as follows:
$$
h(M[r],d)=h(M,d+r),\quad d,r\in\Z,
$$
$$
\HH(M[r],t)=t^{-r}\HH(M,t).
$$
The module $M$ has a minimal finite free graded resolution
\begin{equation}\label{eq1}
0\to \bigoplus_{j\in \Z} R[-j]^{\beta_{r,j}}\to\bigoplus_{j\in \Z}
R[-j]^{\beta_{1,j}}\to\cdots\to
\bigoplus_{j\in \Z}R[-j]^{\beta_{0,j}}\to M\to 0.
\end{equation}
by Hilbert's Syzygy Theorem (see for example \cite{Hilbert},\cite[\S 13]{ZSII},\cite[\S 19.2]{Eisenbud}).
The nonnegative integers $\beta_{i,j}$ are
called the {\it graded Betti numbers}. For all but finitely many pairs $(i,j)$
we have $\beta_{i,j}=0$.
The {\it Betti
numbers\/} are defined by $\beta_i=\sum_{j\in \Z}
\beta_{i,j}$ (not to be confused with the topological Betti numbers
of the complement of the subspace arrangement mentioned earlier).
 Without loss of generality we may assume that $\beta_r\neq 0$.
The nonnegative integer $\cd(M):=r$
is  the {\it cohomological dimension\/} of the module $M$
and is at most $n$. 
The {\it Castelnuovo-Mumford regularity\/} of $M$ is
$$
\reg(M):=\max\{j-i\mid 0\leq i\leq r,\beta_{i,j}\neq 0\}.
$$
From the exactness of (\ref{eq1})
follows that 
\begin{equation}\label{eq11}
h(M,d)=
\sum_{i=0}^r(-1)^i\sum_{j\in \Z}\beta_{i,j}h(R[-j],d)=
\sum_{i=0}^r(-1)^i\sum_{j\in \Z}\beta_{i,j}h(R,d-j)
\end{equation}
and
$$
\HH(M,t)=\frac{\sum_{i=0}^r (-1)^i\sum_{j\in \Z}\beta_{i,j}t^{j}}{(1-t)^n}.
$$
We define the {\it Hilbert polynomial\/}  $\widetilde{h}_M$ of M by
\begin{equation}\label{eq12}
\widetilde{h}_M(d)=\widetilde{h}(M,d)=\sum_{i=0}^r(-1)^i\sum_{j\in \Z}
\beta_{i,j}\widetilde{h}(R,d-j).
\end{equation}
\begin{corollary}\label{cor1}
If
$d\geq
1-n+\reg(M)+\cd(M)$, 
then we have
$$
h_M(d)=\widetilde{h}_M(d).
$$
\end{corollary}
\begin{proof}
If $\beta_{i,j}\neq 0$, then
$$
d-j=d-i-(j-i)\geq d-\cd(M)-\reg(M)\geq 1-n, 
$$
hence
$h_R(d-j)=\widetilde{h}_R(d-j)$
(see Example~\ref{ex1}). 
The corollary follows from (\ref{eq11}) and (\ref{eq12}).
\end{proof}
\section{Subspace arrangements}
For the remainder of this paper, let $V$ be an $n$-dimensional
vector space and suppose that $V_1,V_2,\dots,V_m$
are subspaces of $V$. 
For a subset $X\subseteq V$, let ${\mathcal I}(X)\subseteq R=K[V]$
be its vanishing ideal.
Define $I_j={\mathcal I}(V_j)$ for $j=1,2,\dots,m$.
The union
$$
{\mathcal A}=V_1\cup V_2\cup \cdots\cup V_m
$$
is a subspace arrangement. Its vanishing ideal is
$$
I:={\mathcal I}({\mathcal A})=I_1\cap I_2\cap \cdots \cap I_m.
$$
Define
$$
J:=I_1I_2\cdots I_m.
$$
\begin{theorem}\label{theoregularity}
The regularity of $I$ and $J$ are bounded by
$\reg(I)\leq m$ and
$\reg(J)\leq m$.
\end{theorem}
\noindent For $\reg(J)\leq m$, see~\cite{CH}.
The bound $\reg(I)\leq m$ was proven in \cite{DS1,DS2}.

For a ideal ${\mathfrak a}\subseteq R$ we have
$\cd({\mathfrak a})=\cd(R/{\mathfrak a})-1\leq n-1$.
In particular, we get $\cd(I)\leq n-1$ and $\cd(J)\leq n-1$.
\begin{corollary}\label{cor5}
We have $h_I(d)=\widetilde{h}_I(d)$
and
$h_J(d)= \widetilde{h}_J(d)$ 
for $d\geq m$.
\end{corollary}
\begin{proof}
By Corollary~\ref{cor1}, $h_{I}(d)=\widetilde{h}_{I}(d)$ and
$h_J(d)=\widetilde{h}_J(d)$
for 
$$
d\geq m= (1-n)+m+(n-1)\geq 1-n+\reg(I)+\cd(I).
$$
\end{proof}
\section{Main results}
For $S\subseteq\{1,2,\dots,m\}$, define
$I_S=\bigcap_{s\in S}I_s$ and 
$J_S=\prod_{s\in S}I_s$.
Note that
$I=I_{\{1,2,\dots,m\}}$ and $J=J_{\{1,2,\dots,m\}}$.
We use the convention $I_\emptyset=J_{\emptyset}=R$.
For $S\subseteq\{1,2,\dots,m\}$ define
$V_S=\bigcap_{i\in S}V_i$, $n_S=\dim V_S$ and $c_S=n-n_S$
is the codimension of $V_S$ in $V$. We also set $n_i=n_{\{i\}}=\dim V_i$
and $c_i=c_{\{i\}}=n-n_i$ for $i=1,2,\dots,m$.

We define polynomials $p_S(t)$ recursively as follows. First we define
$$
p_\emptyset(t)=1.
$$
If $S\neq\emptyset$ and $p_X(t)$ is already defined for all proper
subsets $X\subset S$, then 
$p_S(t)$ is uniquely determined by 
$$
\sum_{X\subseteq S}(-t)^{|X|}p_X(t) \equiv 0 \bmod (1-t)^{c_S},\quad
\deg(p_S(t))<c_S.
$$
Here $\deg(p_X(t))$ is the degree of the polynomial $p_X(t)$
and $|X|$ is the cardinality of the set $X$.
\begin{theorem}\label{thmmain1}
We have
$$
\HH(J,t)=\frac{p(t)t^m}{(1-t)^n},
$$
where $p(t)=p_{\{1,2,\dots,m\}}(t)$.
\end{theorem}
\begin{corollary}
The Hilbert series $\HH(J,t)$ depends only on the numbers
$n_S, S\subseteq \{1,2,\dots,m\}$.
\end{corollary}
\begin{proof} This follows immediately from the recursive formulas
for $p_S(t)$, $S\subseteq \{1,2,\dots,n\}$ and the observations
that $n=n_{\emptyset}$ and $c_S=n-n_S$ for all $S$.
\end{proof}
\begin{proposition}\label{prop43}
Let $\{\beta_{i,j}\}$ be the graded Betti numbers of the ideal $J$. If $\beta_{i,j}\neq 0$
then $j=m+i$ and $\beta_i:=\sum_{j}\beta_{i,j}=\beta_{i,m+i}$ is the $i$-th Betti
number.
We have
$$
\HH(J,t)=\frac{\sum_{i=0}^r(-1)^i\beta_it^{i+m}}{(1-t)^n}=\frac{t^mp(t)}{(1-t)^n},
$$
where
$$
p(t)=\beta_0-\beta_1t+\beta_2t^2-\cdots+(-1)^{r}\beta_{r}t^{r},
$$
and $r=\cd(J)\leq n-1$.
The Betti numbers and graded Betti numbers of $J$ are uniquely
determined by the Hilbert series of $J$, so they are combinatorial invariants.
\end{proposition} 
\begin{proof}
By Theorem~\ref{theoregularity} we have
$$
\reg(J)\leq m.
$$
The ideal $J=I_1I_2\cdots I_m$ is generated in degree $m$.
This means that $J$ has a {\it linear\/} minimal free resolution
(see~\cite[Proposition]{EisenbudGoto}):
$$
0\to R[-m-r]^{\beta_r}\to \cdots \to R[-m-1]^{\beta_1}\to R[-m]^{\beta_0}\to J\to 0.
$$
The proposition follows.
\end{proof}
\begin{corollary}\label{corComb}
The Betti numbers and graded Betti numbers of $J$ are uniquely determined
by the Hilbert series of $J$, so they are combinatorial invariants.
\end{corollary}
\begin{definition}
The subspaces $V_1,\dots,V_m$ are called {\it transversal\/}
if
$$
c_S=\min\big(n,\textstyle\sum_{i\in S}c_i\big)
$$
for all $S\subseteq\{1,2,\dots,m\}$, where $\min$ denotes the minimum.
\end{definition}
Note that we always have $c_S\leq \min\big(n,\sum_{i\in S}c_i\big)$.
So the subspaces are transversal if any intersection of some of
the subspaces has the smallest possible dimension.
\begin{theorem}\label{theo1}
Suppose that $V_1,\dots,V_m$ are transversal.
Then $\HH(I,t)-f(t)$ and $\HH(J,t)-f(t)$ are polynomials in $t$,
where
$$
f(t)=\frac{\prod_{i=1}^d \big(1-(1-t)^{c_i}\big)}{(1-t)^n}.
$$
\end{theorem}
\begin{corollary}
If $V_1,\dots,V_m$ are transversal, then
$$
h_I(d)=\widetilde{h}_I(d)=\widetilde{h}_J(d)=h_J(d)
$$
for all $d\geq m$.
\end{corollary}
\begin{proof}
From Theorem~\ref{theo1} and  the exact sequence
$$
0\to J\to I\to I/J\to 0
$$
follows that
 $\HH(I/J,t)=\HH(I,t)-\HH(J,t)$ is a polynomial. So
$I/J$ has a finite dimension. 
It follows that $\widetilde{h}_I(t)=\widetilde{h}_J(t)$.
We have
$$
h_I(d)=\widetilde{h}_I(d)=\widetilde{h}_J(d)=h_J(d)
$$
for $d\geq (1-n)+m+(n-1)=m$ by Corollary~\ref{cor5}.
\end{proof}
\begin{corollary}
If $V_1,\dots,V_m$ are transversal, then
$$
h_I(d)=h_J(d)=
\sum_{S}(-1)^{|S|}
{d+n-1-c_S\choose n-1-c_S}.
$$
where $c_S=\sum_{i\in S}c_i$ and the sum is over all $S\subseteq\{1,2,\dots,m\}$
for which $c_S<n$.
\end{corollary}

\section{Examples}
\begin{example}\label{ex1}
Suppose that $m=3$ and $n=3$ and that $V_1,V_2,V_3$ are one-dimensional
and distinct. Now $V_1,V_2,V_3$ correspond to 3 points 
$P_1,P_2,P_3\in \PP^2$.
Suppose that $P_1,P_2,P_3$ are not colinear. We can change coordinates
such that 
$$P_1=(1:0:0),\ P_2=(0:1:0),\ P_3=(0:0:1).
$$
We have
$$
I=I_1\cap I_2\cap I_3=(x_2,x_3)\cap (x_1,x_3)\cap (x_1,x_2)=
(x_2x_3,x_1x_3,x_1x_2)\subseteq K[x_1,x_2,x_3]
$$
and
$$
J=I_1I_2I_3=(x_1x_2x_3,x_1^2x_2,x_1^2x_3,x_2^2x_3,x_2^2x_1,x_3^2x_1,x_3^2x_2)\subseteq
K[x_1,x_2,x_3].
$$
The minimal free resolutions are as follows
$$
0\to R[-3]^2\to R[-2]^3\to I\to 0,
$$
$$
0\to  R[-5]^3\to R[-4]^9\to R[-3]^7 \to J\to 0.
$$
So we have $\cd(I)=1$, $\cd(J)=2$, $\reg(I)=2$, $\reg(J)=3$,
$$
\HH(I,t)=\frac{3t^2-2t^3}{(1-t)^3}\mbox{ and }
\HH(J,t)=\frac{7t^3-9t^4+3t^5}{(1-t)^3}.
$$
The Hilbert functions, which can be found as the
coefficients of the Hilbert series,  are given in the following table.
$$
\begin{array}{c||c|c|c|c|c|c}
d & 0 & 1 & 2 & 3 & 4 & 5 \\ \hline\hline
h_{I}(d) & 0  & 0 & 3 & 7 & 12 & 18\\ \hline
h_J(d) & 0 & 0 & 0 & 7 & 12 & 18
\end{array}
$$
The Hilbert polynomials $\widetilde{h}_I(d)$ and $\widetilde{h}_J(d)$
are both equal to
$$
\frac{d^2+3d-4}{2}.
$$

Let us verify Theorems~\ref{thmmain1} and \ref{theo1}. The subspaces $V_1,V_2,V_3$
are transversal. Let 
\begin{equation}\label{eq55}
f(t)=\prod_{i=1}^3\frac{(1-(1-t)^{c_i})}{(1-t)^n}=\frac{(1-(1-t)^2)^3}{(1-t)^3}=
\frac{-2+6t-3t^2}{(1-t)^3}+2-3t^2+t^3.
\end{equation}
We can write
$$
\HH(I,t)=\frac{3t^2-2t^3}{(1-t)^3}=\frac{-2+6t-3t^2}{(1-t)^2}+2
$$
and
$$
\HH(J,t)=\frac{7t^3-9t^4+3t^5}{(1-t)^3}=\frac{-2+6t-3t^2}{(1-t)^2}+2-3t^2.
$$
As Theorem~\ref{theo1} predicts, $\HH(I,t)-f(t)$ and $\HH(J,t)-f(t)$ are polynomials.

Using the recursive definitions of $p_S$ we get $p_{\emptyset}(t)=1$, 
$$
p_{\emptyset}(t)-tp_{\{1\}}(t)\equiv 1-t p_{\{2\}}(t)\equiv 1-tp_{\{3\}}(t)\equiv 0 \bmod (1-t)^2,
\quad \deg(p_{\{1\}}),\deg(p_{\{2\}}),\deg(p_{\{3\}})<2,
$$
so $p_{\{1\}}(t)=p_{\{2\}}(t)=p_{\{3\}}(t)=2-t$. Furthermore, 
$$
1-t(p_{\{1\}}(t)+p_{\{2\}}(t))+t^2p_{\{1,2\}}(t)\equiv 0\bmod (1-t)^3,
$$
so we have
$$
t^2p_{\{1,2\}}(t)\equiv -1+2t(2-t)\equiv -1+4t-2t^2 \bmod (1-t)^3
$$
and $\deg(p_{\{1,2\}})<3$.
It follows that
$$
p_{\{1,2\}}(t)=p_{\{1,3\}}(t)=p_{\{2,3\}}(t)=4-4t+t^2
$$
Finally we have
$$
1-t(p_{\{1\}}(t)+p_{\{2\}}(t)+p_{\{3\}}(t))+
t^2(p_{\{1,2\}}(t)+p_{\{1,3\}}(t)+p_{\{2,3\}}(t))-t^3p_{\{1,2,3\}}(t)\equiv
0 \bmod (1-t)^3,
$$
so
$$
t^3p_{\{1,2,3\}}(t)\equiv 1-3t(2-t)+3t^2(4-4t+t^2)\equiv  -2+6t-3t^2\bmod (1-t)^3.
$$
From this follows that
$$
p_{\{1,2,3\}}(t)=7-9t+3t^2.
$$
Theorem~\ref{thmmain1} correctly gives
$$
\HH(J,t)=\frac{t^3p_{\{1,2,3\}}(t)}{(1-t)^3}=
\frac{7t^3-9t^4+3t^5}{(1-t)^3}.
$$
\end{example}

\begin{example}\label{ex2}
If $P_1,P_2,P_3$ are colinear, then, after a change of coordinates, 
we may assume
$$
P_1=(1:0,0),\ P_2=(0:1:0),\ P_3=(1:1:0).
$$
We have
$$
I=I_1\cap I_2\cap I_2=(x_2,x_3)\cap (x_1,x_3)\cap (x_1-x_2,x_3)=
(x_1x_2(x_1-x_2),x_3)\subseteq K[x_1,x_2,x_3]
$$
and
$$
J=I_1I_2I_3=(x_1x_2(x_1-x_2),x_1^2x_3,x_1x_2x_3,x_2^2x_3,x_1x_3^2,x_2x_3^2,x_3^3)\subseteq
K[x_1,x_2,x_3].
$$
We have minimal free resolutions
$$
0\to R[-4]\to R[-1]\oplus R[-3]\to I\to 0,$$
$$
0\to R[-5]^3\to R[-4]^6\to R[-3]^7\to J\to 0.
$$
We get $\cd(I)=1$, $\cd(J)=2$, $\reg(I)=\reg(J)=3$,
$$
\HH(I,t)=\frac{t+t^3-t^4}{(1-t)^3}\mbox{ and }
\HH(J,t)=\frac{7t^3-9t^4+3t^5}{(1-t)^3}.
$$
The Hilbert functions are
$$
\begin{array}{c||c|c|c|c|c|c}
d & 0 & 1 & 2 & 3 & 4 & 5 \\ \hline\hline
h_{I}(d) & 0  & 1 & 3 & 7 & 12 & 18\\ \hline
h_J(d) & 0 & 0 & 0 & 7 & 12 & 18
\end{array}
$$
We can compute $f(t)$ and it is the same as in Example~\ref{ex1}, (\ref{eq55}).
Now
$$
\HH(I,t)=\frac{t+t^3-t^4}{(1-t)^3}=\frac{-2+6t-3t^2}{(1-t)^3}+2+t,
$$
so $\HH(I,t)-f(t)$ is indeed a polynomial as in Theorem~\ref{theo1}.
Theorem~\ref{thmmain1} gives the same result for $\HH(J,t)$ as in
Example~\ref{eq1}.
\end{example}
%
\begin{example}\label{ex3}
Let $n=4$ and $m=3$.
Consider the ideals $I_1=(x_2,x_3)$, $I_2=(x_1,x_3)$ and $I_3=(x_1,x_2)$
in the polynomial ring $K[x_1,x_2,x_3,x_4]$ in 4 variables
and let again $I=I_1\cap I_2\cap I_3$.
The ideals $I_1,I_2,I_3$ correspond to 2-dimensional
subspaces $V_1,V_2,V_3$ of the 4-dimensional space $V$ such that 
$V_1\cap V_2\cap V_3$ is $1$-dimensional, and $V_1+V_2+V_3=V$.
We get similar free resolutions of $I$ as in Example~\ref{ex1}:
$$
0\to R[-3]^2\to R[-2]^3\to I\to 0,
$$
We have
$$
\HH(I,t)=\frac{3t^2-2t^3}{(1-t)^4}.
$$
\end{example}
\begin{example}\label{exlast}
Let $n=4$ and $m=3$.
Consider the ideals $I_1=(x_1,x_3)$, $I_2=(x_2,x_3)$ and $I_3=(x_1-x_2,x_3)$
in the polynomial ring $K[x_1,x_2,x_3,x_4]$ in 4 variables
and let $I=I_1\cap I_2\cap I_3$, $J=I_1I_2I_3$.
The ideals $I_1,I_2,I_3$ correspond to 2-dimensional
subspaces $V_1,V_2,V_3$ of the 4-dimensional space $V$ such that 
$V_1\cap V_2\cap V_3$ is $1$-dimensional, and $V_1+V_2+V_3$ is 3-dimensional.
In projective space $\PP^3$, we have 3 lines lying an a plane
and going through 1 point.
We get similar free resolutions of $I$ and $J$ as in Example~\ref{ex2}.
$$
0\to R[-4]\to R[-1]\oplus R[-3]\to I\to 0,
$$
$$
0\to R[-5]^3\to R[-4]^9\to R[-3]^7\to J\to 0.
$$
We obtain
$$
\HH(I,t)=\frac{t+t^3-t^4}{(1-t)^4}\mbox{ and
}\HH(J,t)=\frac{7t^3-9t^4+3t^5}{(1-t)^4}.
$$
Taking the difference gives
$$
\HH(I,t)-\HH(J,t)=\frac{3t^2+t}{1-t},
$$
which is not a polynomial. Note that Theorem~\ref{theo1} does not apply
because $V_1,V_2,V_3$ are not transversal.
We have
$$
\HH(I,t)=\frac{t+t^3-t^4}{(1-t)^4},
$$
but in Example~\ref{ex3} we got
$$
\HH(I,t)=\frac{3t^2-2t^3}{(1-t)^4}.
$$
This shows that $\HH(I,t)$ is not a combinatorial invariant.
The difference of both Hilbert series is not even a polynomial.
This implies that the Hilbert {\it polynomial\/} $h_I(d)$ is not
a combinatorial invariant either.
\end{example}
\section{Complexes of product ideals and intersection ideals}
\begin{theorem}[See Chapter IV of \cite{Sidman}]\label{theoSidman}
There exists complexes
$$
0\to I\to \bigoplus_{|S|=m-1}I_{S}\to
\bigoplus_{|S|=m-2}I_S\to \cdots \to\bigoplus_{|S|=1}I_S\to R\to 0
$$
and
$$
0\to J\to \bigoplus_{|S|=m-1}J_{S}\to
\bigoplus_{|S|=m-2}J_S\to \cdots \to \bigoplus_{|S|=1}J_S\to R\to 0.
$$
whose homologies are killed by ${\mathfrak a}=\sum_{j=1}^m I_j$.
\end{theorem}
To describe the the maps in the complexes in Theorem~\ref{theoSidman}
it suffices to define maps $I_T\to I_S$ and $J_T\to J_S$
for all subsets $S,T\subseteq \{1,2,\dots,m\}$ with $|T|=|S|+1$.
If $T=\{i_1,i_2,\dots,i_r\}$ with $i_1<i_2<\cdots<i_r$
and $S=\{i_1,i_2,\dots,i_{s-1},i_{s+1},\dots,i_r\}$
then the maps $I_T\to I_S$ and $J_T\to J_S$ 
in the complexes in Theorem~\ref{theoSidman} are given
by $f\mapsto (-1)^sf$. All other maps are equal to $0$.

\begin{corollary}\label{corHilbertseries}
If $V_{\{1,2,\dots,m\}}=\bigcap_{i=1}^m V_i=(0)$, then
$$
\sum_{S\subseteq \{1,2,\dots,m\}}(-1)^{|S|}\HH(I_S,t)
$$
and
$$
\sum_{S\subseteq \{1,2,\dots,m\}}(-1)^{|S|}\HH(J_S,t)
$$
are polynomials in $t$.
\end{corollary}
\begin{proof}
The ideal
$$
\sum_{j=1}^m I_j={\mathcal I}\big({\textstyle \bigcap_{i=1}^m V_i}\big)=
{\mathcal I}(\{0\})={\mathfrak m}
$$
is the maximal homogeneous ideal of $R$.

Suppose that 
$$
\xymatrix{
0\ar[r]^{\partial_{r+1}} & C_{r}\ar[r]^{\partial_r} &
C_{r-1}\ar[r]^{\partial_{r-1}} &\cdots\ar[r]^{\partial_1} & 
C_{0}\ar[r]^{\partial_0} & 0}
$$
is a complex of finitely generated graded $R$-modules. 
The $i$-th homology group is
$$H_i=\ker(\partial_{i})/\im(\partial_{i+1}).
$$ 
We have
$$
\sum_{i=0}^r (-1)^i \HH(C_i,t)=\sum_{i=0}^r(-1)^i \HH(H_i,t).
$$
If  ${\mathfrak m}H_i=0$, then $H_i$ is
finite dimensional, and $\HH(H_i,t)$ is a polynomial for all~$i$.

We now apply this to the complexes in Theorem~\ref{theoSidman}.
\end{proof}
\begin{proof}[Proof of Theorem~\ref{thmmain1}]
By Proposition~\ref{prop43} we can write
$$
\HH(J,t)=\frac{\sum_{i=0}^r(-1)^i\beta_it^{i+m}}{(1-t)^n}=\frac{t^mp(t)}{(1-t)^n},
$$
where
\begin{equation}\label{eqBetti}
p(t)=\beta_0-\beta_1t+\cdots+(-1)^r\beta_rt^r
\end{equation}
is a polynomial of degree $r\leq \cd(J)\leq n-1$. Similarly we
can write
$$
\HH(J_S,t)=\frac{t^{|S|}p_S(t)}{(1-t)^n}
$$
with 
\begin{equation}\label{eqdegbound}
\deg(p_S(t))\leq n-1
\end{equation}
for all $S\subseteq \{1,2,\dots,m\}$.

Let $W=V/V_{S}$ and define $W_X=V_X/V_{S}$ for $X\subseteq S$.
Let $\overline{J}_i\subseteq K[W]\cong K^{[n-n_S]}$ be the vanishing ideal
of $W_i$ for all $i\in S$.  Define $\overline{J}_X=\prod_{i\in
X}\overline{J}_i$ for all $X\subseteq S$.

We have
$$
\overline{J}_X\otimes K^{[n_S]}=J_X
$$
inside $K^{[n]}=K^{[n-n_S]}\otimes K^{[n_S]}$.
From this follows that
$$
\frac{t^{|X|}p_X(t)}{(1-t)^n}=
\HH(J_X,t)=\frac{\HH(\overline{J}_X,t)}{(1-t)^{n_S}}.
$$
In particular, we have
$$
\HH(\overline{J}_S,t)=\frac{t^{|S|}p_S(t)}{(1-t)^{n-n_S}}.
$$
From this it follows that $\deg(p_S(t))\leq \dim W-1=n-n_S-1$
(see~(\ref{eqdegbound})).
Since $\bigcap_{i\in S}W_i=0$ in $W$, Corollary~\ref{corHilbertseries} implies
that
$$
\sum_{X\subseteq S}(-1)^{|X|}\HH(\overline{J}_X,t)=\sum_{X\subseteq S}
\frac{(-t)^{|X|}p_X(t)}{(1-t)^{n-n_S}}
$$
is a polynomial in $t$.
Multiplying with $(1-t)^{n-n_S}$ gives 
$$
\sum_{X\subseteq S}(-t)^{|X|}p_X(t)\equiv 0 \bmod (1-t)^{n-n_S}.
$$
\end{proof}

\begin{proof}[Proof of Theorem~\ref{theo1}]\

{\bf  Special case:} Suppose that $c_1+\cdots+c_m<n$. 
After a change of coordinates (and using that the
arrangement is transversal) we may identify
$K^{[n]}$
with
$$
K^{[c_1]}\otimes K^{[c_2]}\otimes \cdots \otimes K^{[c_m]}\otimes
K^{[n-c_1-\cdots-c_m]}
$$
and
$I_k$ with
$$
K^{[c_1]}\otimes\cdots\otimes K^{[c_{k-1}]}\otimes {\mathfrak m}_{c_k}\otimes 
K^{[c_{k+1}]}\otimes \cdots\otimes K^{[c_{m}]}\otimes
K^{[n-c_1-\cdots-c_m]}
$$
for
$k=1,2,\dots,m$.
Here ${\mathfrak m}_r$ is the homogeneous maximal ideal of $K^{[r]}$.
We get
$$
I=
{\mathfrak m}_{c_1}\otimes {\mathfrak m}_{c_2}\otimes \cdots\otimes
{\mathfrak m}_{c_m}\otimes K^{[n-c_1-\cdots-c_m]}=J.
$$
Note that $\HH(K^{[r]},t)=(1-t)^{-r}$ and $\HH({\mathfrak m}_r,t)=(1-t)^{-r}-1$ 
for
all $r$.
Therefore, we get
$$
\HH(I,t)=
\HH(J,t)=H({\mathfrak m}_{c_1},t)\HH({\mathfrak m}_{c_2},t)\cdots 
\HH({\mathfrak m}_{c_m},t)\HH(K^{[n-c_1-\cdots-c_m]},t)=
$$
$$
=\left(\prod_{i=1}^m \big((1-t)^{-c_i}-1\big)\right) \cdot (1-t)^{c_1+\cdots+c_m-n}=
\frac{\prod_{i=1}^m \big(1-(1-t)^{c_i})}{(1-t)^n}.
$$

\noindent {\bf The general case.} We prove the Theorem by induction on $m$.
The base case $m=1$ follows from the special case above. If $c_1+\cdots+c_m<n$
then we are also done. Let us assume that $c_1+\cdots+c_m \geq n$.
In particular we have $c_{\{1,2,\dots,m\}}=\min(n,\sum_{i=1}^mc_i)=n$, so
$\bigcap_{i=1}^mV_i=(0)$. By Corollary~\ref{corHilbertseries} we have that
$$
\sum_{S\subseteq \{1,2,\dots,m\}}(-1)^{|S|}\HH(I_S,t)
$$
is a polynomial.

By induction we have that
$$\HH_{I_S}(t)-w\prod_{i\in S}v_i
$$
is a polynomial for all strict subsets $S\subset \{1,2,\dots,m\}$,
where $v_i=1-(1-t)^{c_i}$ and  $w=(1-t)^{-n}$.
To show that
$$
\HH_I(t)-w\prod_{i=1}^mv_i
$$
is a polynomial, it suffices to show that
\begin{equation}\label{eqbigsum}
\sum_{S\subseteq \{1,2,\dots,m\}}(-1)^{|S|}\big(\HH(I_S,t)-w\prod_{i\in
S}v_i\big)
\end{equation}
is a polynomial.

Now
$$
\sum_{S\subseteq\{1,2,\dots,m\}}(-1)^{|S|}\HH(I_S,t)
$$
is a polynomial by Corollary~\ref{corHilbertseries}, and
$$
\sum_{S\subseteq\{1,2,\dots,m\}}(-1)^{|S|}w\prod_{i\in  S}v_i=
w\prod_{i=1}^m(1-v_i)=(1-t)^{-n}\prod_{i=1}^m(1-t)^{c_i}=(1-t)^{c_1+\cdots+c_m-n}
$$
is a polynomial as well.
Therefore (\ref{eqbigsum}) is a polynomial.
\end{proof}
\section{Application to Generalized Pricipal Component Analysis}
The object of Principal Component Analysis (PCA) is to
approximate a data set inside a vector space $V$ by a subspace
of smaller dimension. In Generalized Principal Component
Analysis (GPCA) one tries to approximate a data set inside a vector space
$V$ by a {\it union\/} of subspaces spaces (in other words, a {\it subspace arrangement}). 
Some applications of GPCA are
 motion segmentation (see~\cite{VidalR2002,VidalR2004}),
  image segmentation (see~\cite{VidalR2004E}), 
  image compression (see~\cite{HongW}) and hybrid control
  systems (\cite{MaY}). For an overview of GPCA, see~\cite{VMS}.

A first start in GPCA is to decide on the number of subspaces
and the dimensions of the subspaces of the subspaces arrangement
that will approximate the data.

Suppose that $v_1,v_2,\dots,v_r\in V$ are data points. 
Here $r$ is fairly large.
Suppose that $v_1,\dots,v_r$ are contained in some
subspace arrangement ${\mathcal A}=V_1\cup \cdots \cup V_m$,
unknown to us.
We would like to recover $n_1,\dots,n_m$ where $n_i=\dim V_i$.
Let ${\mathfrak a}_j$ be the vanishing ideal of the ray through $v_j$.
Then we have that
$$
h({\mathfrak a}_1\cap \cdots\cap {\mathfrak a}_r,d)=
h(I,d)
$$
for small values of $d$, where $I={\mathcal I}({\mathcal A})$ as before.
Now 
$$
h({\mathfrak a}_1\cap \cdots\cap {\mathfrak a}_r,d)
$$
can be computed using linear algebra for small values of $d$.
Therefore, we can determine $h(I,d)$ for small values
of $d$. So an important question is,
given $h(I,d)$ for small values of $d$,
can we determine the dimensions $n_1,n_2,\dots,n_m$?
Proposition~\ref{prop18} gives an affirmative answer if the
subspaces are transversal. Of course, in real applications
the data is approximated by the subspaces arrangement, but
not contained in it. In that case, using the PCA method
in $K[V]_d$ one can still can estimate the value
$h({\mathcal I}({\mathcal A}),d)$.
\begin{proposition}\label{prop18}
Assume that the arrangement is transversal.
Suppose that $c_1,\dots,c_m$ are unknown, but we know the values of the
Hilbert polynomial
$$
h_I(d)
$$
for $d=m,m+1,\dots,m+n-1$, then we can recover $c_1,\dots,c_m$.
\end{proposition}
\begin{proof}
Note that $h_I(d)=\widetilde{h}_I(d)$ for $d\geq m$.
Since we know $\widetilde{h}_I(d)$ for $d=m,m+1,\dots,m+n-1$
and $\widetilde{h}_I$ has degree $\leq n-1$, $\widetilde{h}_I$ is
uniquely determined. From this, we can determine $\HH(I,t)$, up to a polynomial.
Suppose that 
$\HH(I,t)$
is equal to $a(t)/(1-t)^n$ up to a polynomial.
Let $b(t)$ be the remainder of division of $a(t)$ by $(1-t)^n$.
Then $b(t)$ has degree $<n$ and $\HH(I,t)$ is equal to $b(t)/(1-t)^n$
modulo a polynomial. So we have that
$$
b(t)\equiv \prod_{i=1}^d(1-(1-t)^{c_i})\bmod (1-t)^n 
$$
and 
$$
b(1-t)\equiv \prod_{i=1}^d(1-t^{c_i})\bmod t^n.
$$
Let $r_i$ be the number of the $c_j$'s equal to $i$.
Then we have
$$
b(1-t)\equiv \prod_{i=1}^d(1-t^i)^{r_i}\bmod t^n.
$$
From this we can easily determine $r_1,r_2,\dots,r_{n-1}$ (in that order).
Indeed, if we already know $r_1,\dots,r_s$, then the Taylor series
of
$$
\frac{b(1-t)}{\prod_{i=1}^{s}(1-t^i)^{r_i}}
$$
is
$$
1-r_{s+1}t^{s+1}+\mbox{higher order terms.}
$$
So we find the value of $r_{s+1}$.
\end{proof}

\end{document}